\documentclass[11pt,a4paper,reqno]{article}
\usepackage[utf8]{inputenc}
\usepackage{mathptmx}
\usepackage[numbers]{natbib}

\usepackage{amsfonts,amsmath,amssymb,amsthm}       
\usepackage{enumerate}

\numberwithin{equation}{section}

\newtheorem{theorem}{Theorem}[section]
\newtheorem{corollary}[theorem]{Corollary}
\newtheorem{proposition}[theorem]{Proposition}
\newtheorem{lemma}[theorem]{Lemma}

\theoremstyle{definition}
\newtheorem{example}[theorem]{Example}
\newtheorem{definition}[theorem]{Definition}
\newtheorem{remark}[theorem]{Remark}

\author{Jaakko Lehtomaa\thanks{Department of Mathematics and Statistics, University of Helsinki, 
P.O. Box 68, 00014 Helsinki, Finland, 
e-mail address: jaakko.lehtomaa@helsinki.fi }}

\title{A note on limiting behaviour of constrained sums of two variables}

\begin{document}
\maketitle

\begin{abstract}
This note studies the asymptotic properties of the variable 
$$Z_d:=\frac{X_1}{d}\big{|}\{X_1+X_2=d\},$$
as $d\to \infty$. Here $X_1$ and $X_2$ are non-negative i.i.d. variables with a common twice differentiable density function $f$.

General results concerning the distributional limits of $Z_d$ are discussed with various examples. Eventual log-convexity or log-concavity of $f$ turns out to be the key ingredient that determines how the variable $Z_d$ behaves. As a consequence, two surprising discoveries are presented: Firstly, it is noted that the distributional limit is not strictly determined by the decay rate of the tail function. Secondly, it is shown that there exists a light-tailed distribution exhibiting behaviour that is commonly associated with heavy-tailed distributions i.e. the principle of a single big jump.
\end{abstract}

\noindent \emph{MSC classification (2010)}: 60E05; 60F05; 62E20


\section{Preliminaries}

During the last decades it has become clear that heavy-tailed random variables are needed in realistic mathematical models. Consequently, heavy-tailed analysis has seen an explosive growth in the number of publications, making it an active research field of high current interest. 

A cornerstone of heavy-tailed thinking is the \emph{principle of a single big jump}. Unfortunately, there does not seem to exist consensus about the exact definition of this principle. Nevertheless, the principle always consists of the idea that the most likely way for a sum to be large is that one of the summands is large. Some authors refer to this principle whenever there exists a dominating random variable \cite{embrechts:1997,foss:2007}, whereas other reserve the expression for subexponential distributions \cite{armendariz:2011,asmussen:1996,denisov:2008} or their generalisations \cite{beck:2015}. Some properties are also studied in the case of dependent variables \cite{albrecher:2006}. 

The aim of this note is to study the principle of a single big jump in a rigorous setting. In \cite{foss:2013}, the behaviour of the process $(Z_d):=(Z_d)_{d>0}$ is used to illustrate the phenomenon of a single big jump. Our plan is to study the process $(Z_d)$ further and to present general results whose applicability can be verified using the density function $f$. 

In order to do this, we define two convergence types for the process $(Z_d)$:
\begin{enumerate}[I)]
\item $\mathcal{L}(Z_{d})\to\frac{1}{2} \delta_0+\frac{1}{2} \delta_1$ \label{a} and
\item $\mathcal{L}(Z_{d})\to \delta_{\frac{1}{2}}$.\label{b}
\end{enumerate}
In \ref{a} and \ref{b} the notation $\mathcal{L}(Z_{d})$ refers to the law of $Z_d$ and the convergence is understood as convergence in distribution in the limit $d\to \infty$. In Types \ref{a} and \ref{b}, $\delta_x$ signifies a distribution concentrated to the point $x\in \{0,1/2,1\}$. 

Behaviour \ref{a} resembles the way many heavy-tailed variables are known to behave: if the sum $X_1+X_2$ is large then one of the variables is large. Behaviour \ref{b} is related to a phenomenon encountered within the class of light-tailed distributions: both of the variables $X_1$ and $X_2$ contribute equally. 

Recall that a random variable $X$ is called \emph{heavy-tailed} if $E(e^{sX})=\infty$ for all $s>0$ and light-tailed otherwise. We will show that, in the sense of Behaviour \ref{a}, the principle can occur outside the class of heavy-tailed distributions. Traditionally the idea of the principle of a single big jump is almost exclusively associated with a subclass of heavy-tailed distributions called subexponential distributions. The subexponential class and its extensions are further discussed in Section \ref{discussion} below.

\subsection{Assumptions}\label{assumptions}

The non-negative random variables $X_1$ and $X_2$ are independent and identically distributed. The variable $X_1$ has an unbounded support and a density function $f$. Set $F(x):=P(X_1\leq x)$ and $\overline{F}(x):=1-F(x)$. The function $f$ is assumed to be twice differentiable in the set $[0,\infty)$ and eventually decreasing. A property is said to hold \emph{eventually} if there exists  $y_0\in \mathbb{R}$ such that the property is valid in the set $[y_0,\infty)$.

\subsection{Basic Properties}

The density function $f_{Z_d}$ of the variable $Z_d$ can be directly obtained from the conditional distribution of $X_1|\{X_1+X_2\}$. Its density is concentrated in the interval $[0,1]$ and given by formula
\begin{equation}\label{f1}
f_{Z_d}(x)=\frac{f(d x)f(d(1-x))}{\displaystyle\int_{0}^1 f(d y)f(d(1-y)) \, dy},\quad x\in [0,1].
\end{equation}
The function $f_{Z_d}$ can be viewed as a function of two variables as 
$$g(x,d):=f_{Z_d}(x)\colon [0,1] \times (0,\infty)\to [0,\infty).$$
For a fixed $d>0$ the function $f_{Z_d}(x)$ is symmetric with respect to the point $x=1/2$. Hence, it suffices to formulate the results only for  $x\in[0,1/2]$.

Conditions implying Behaviours \ref{a} or \ref{b} typically involve estimation of decay rates of integrals. What is more, neither of the behaviours needs to occur; the distributional limit may exist without any concentration of probability mass. To see this, consider the following example.

\begin{example}\label{ex1}
Suppose $f$ is a gamma density function $f(x)=Cx^{a-1}e^{-x}$, where $x>0$, $a>0$ and $C>0$ is an integration constant.

Then $f_{Z_d}$ of \eqref{f1} reduces to
$$f_{Z_d}(x)=\frac{x^{a-1}(1-x)^{a-1}}{\int_0^1 y^{a-1}(1-y)^{a-1} \, dy},$$
for all $d>0$. So, $\mathcal{L}(Z_d)$ does not depend on $d$ and belongs to the family of Beta distributions.
\end{example}

In order to understand the behaviour of the process $(Z_d)$ one needs additional assumptions to those made in Section \ref{assumptions}. One way to proceed is to demand that the function $f_{Z_d}$ should eventually stay convex or concave at the midpoint of $[0,1]$.  This leads to the following characterisation.

\begin{lemma}\label{lchar}Suppose
\begin{equation}\label{L}
L:=\lim_{x \to \infty} \textnormal{sign}\left( \frac{d^2}{dx^2} \log f(x) \right)
\end{equation}
exists, where 
\begin{displaymath}
   \textnormal{sign}(x) := \left\{
     \begin{array}{rl}
       1 & : x>0 \\
       0 & : x =0 \\
       -1 & :x<0.
     \end{array}
   \right.
\end{displaymath} 

Then the function $f_{Z_d}$ of Formula \eqref{f1} is eventually, in $d$, strictly convex with respect to the variable $x$ at point $x=1/2$ if and only if $L=1$. Similarly, $f_{Z_d}$ is eventually, in $d$, strictly concave  with respect to the variable $x$  at point $x=1/2$ if and only if $L=-1$.
\begin{proof}
Consider the eventually convex case; the eventually concave case is analogous. Let $d>0$. For any $x\in (0,1)$, 
\begin{eqnarray*}
f_{Z_d}''(x)&=&\frac{d^2}{\int_0^1 f(dy)f(d(1-y)) \, dy} [ f''(dx)f(d(1-x))-f'(dx)f'(d(1-x)) \\
&-&f'(dx)f'(d(1-x))+f(dx)f''(d(1-x))].
\end{eqnarray*}
The requirement $f_{Z_d}''(1/2)>0$ simplifies to $f''(d/2)f(d/2)-f'(d/2)^2>0. $
This is equivalent with the statement 
\begin{equation}\label{limi}
\left( \frac{d^2}{dx^2} \log f(x)\right)_{|x=d/2}>0.
\end{equation}
The claim follows upon noticing that $L=1$ holds if and only if \eqref{limi} holds eventually in $d$.
\end{proof}
\end{lemma}

\begin{remark} If $L=0$ in \eqref{L}, then $f''(x)f(x)-f'(x)^2=0$ eventually. The function $f(x)=C_1e^{C_2 x}$ solves this differential equation. Here, $C_1,C_2\in \mathbb{R}$ are suitable constants. Direct application of \eqref{f1} shows that $Z_d$ is eventually uniformly distributed.
\end{remark}

\subsection{Relation of Log-convexity and Log-concavity to Failure Rates}

The condition \eqref{L} implies eventual convexity or concavity of $f$.

\begin{definition} A twice differentiable function $g$ is said to be \emph{eventually strictly convex} if there exists a number $x_0>0$  such that $g''(x)>0$ for all $x>x_0$. \emph{Eventually strictly concave} functions are defined similarly.
\end{definition}

If $L=1$ ($L=-1$) in Equation \eqref{L}, the function $f$ is eventually strictly log-convex (log-concave). This is equivalent with the function $f'(x)/f(x)$ being eventually strictly increasing (decreasing). 

Proceeding as in Lemma 4 of \cite{bagnoli:2005} one obtains for eventually strictly log-convex $f$ and for any $x>x_0$ that
\begin{equation}\label{barv1}
\frac{f'(x)}{f(x)}\int_x^\infty f(y) \, dy<\int_x^\infty \frac{f'(y)}{f(y)}f(y) \, dy. 
\end{equation}
Straightforward calculation reveals Equation \eqref{barv1} being equivalent with
\begin{equation}\label{dfr}
\frac{d}{dx} \left( \frac{f(x)}{\overline{F}(x)} \right)>0.
\end{equation}
Equation \eqref{dfr} implies that the \emph{failure rate} $f(x)/\overline{F}(x)$ is an eventually strictly increasing function and that $\overline{F}$ is a an eventually strictly log-convex function. It can be shown similarly that eventually strictly log-concave densities lead to eventually strictly decreasing failure rates and eventual strict log-concavity of the function $\overline{F}$.

The log-concavity and log-convexity are known be the determining properties in several economical, statistical, probabilistic and operations research related concepts. These classical properties are closely linked, as shown above, to the variables whose failure rate is increasing or decreasing. For additional properties the reader is referred to \cite{hansen:1988,wang:1986,lariviere:2006,banciu:2013}. 

In the current note a new phenomenon where log-convexity or log-concavity plays a central role is discovered. It is the deciding factor that determines the eventual shape of the density of $Z_d$.

\section{The Main Result and Applications}

As mentioned earlier, additional conditions need to be imposed in order to obtain Behaviour \ref{a} or \ref{b}. The first result, Proposition \ref{typechar}, does exactly this, but it requires that $f_{Z_d}(x)\to 0$ for all $x\in(0,1/2)$. This may be tedious to check unless the density is extremely simple. However, the latter result, Theorem \ref{ekalause}, provides a sufficient condition which guarantees the validity of the required property.

\subsection{Theoretical Results}\label{tresults}

\begin{proposition}\label{typechar} Suppose the limit $L$ of Equation \eqref{L} exists. Assume further that $f_{Z_d}(x)\to 0$ for all $x\in (0,1/2)$, as $d\to \infty$.

If $L=1$, then \ref{a} holds. If $L=-1$, then \ref{b} holds.
\begin{proof} Suppose $L=-1$. Now, there exists a number $x_0$ such that for all $x>x_0$ it holds that 
\begin{equation}\label{derkaava1}
\frac{d}{dx} \left( \frac{f'(x)}{f(x)} \right)<0.
\end{equation}
Equation \eqref{derkaava1} implies that the function $f'(x)/f(x)$ is strictly decreasing for $x>x_0$. 

Suppose $d>2 x_0$. Direct calculation reveals that $f_{Z_d}'(x)=0$ if and only if
\begin{equation}\label{critical}
\frac{f'(dx)}{f(dx)}=\frac{f'(d(1-x))}{f(d(1-x))}.
\end{equation}
Therefore, the point $x=1/2$ is always a critical point. In addition, based on Equation \eqref{derkaava1} and symmetry, there are no other critical points in the interval $[x_0/d,1/2]$ i.e. the function $f_{Z_d}$ is monotone in the interval $[x_0/d,1/2]$. The critical point at $x=1/2$ must be a maximum, because $f_{Z_d}''(1/2)<0$. Hence the function $f_{Z_d}$ is increasing in the interval $[x_0/d,1/2]$

Next, it is shown that $f_{Z_d}\to 0$ uniformly in the set $[0,c]$, where $c\in (0,1/2)$. It suffices to note the following. Set
$$M:= \frac{1}{f(x_0)}\max_{z \in [0,x_0]} f(z). $$
Recall that $f$ is continuous and eventually decreasing. Let $y_0\in \mathbb{R}$ be chosen so that $f(x)$ is decreasing when $x>y_0$. Now, for any $x\in[0,x_0/d]$ and all $d$ satisfying $d(1-x_0/d)>y_0$ it holds that
\begin{eqnarray}
f_{Z_d}(x)&=&\frac{f(d x)f(d(1-x))}{\int_{0}^1 f(d y)f(d(1-y)) \, dy} \nonumber\\
&\leq& \frac{\left( \max_{z \in [0,x_0]} f(z) \right)f(d(1-x))}{\int_{0}^1 f(d y)f(d(1-y)) \, dy} \nonumber\\
&\leq& \frac{\left( \max_{z \in [0,x_0]} f(z) \right)f(d(1-x_0/d))}{\int_{0}^1 f(d y)f(d(1-y)) \, dy} \nonumber\\
&=& M f_{Z_d}(x_0/d) \nonumber\\
&\leq& M f_{Z_d}(c).\label{vikaper}
\end{eqnarray}
The right hand side of \eqref{vikaper} converges to $0$, as $d\to \infty$. This implies the desired uniform convergence. 

The uniform convergence and symmetry of the function $f_{Z_d}$ with respect to the point $x=1/2$ imply that for any open set $A\subset [0,1]$ one has
$$\liminf_{d \to \infty} P(Z_d\in A)\geq \delta_\frac{1}{2}(A).$$
This is precisely the Portmanteau characterisation of the distributional convergence and the proof of the case $L=-1$ is complete.

If $L=1$, the proof is simpler. In this case the monotonicity together with the assumption $f_{Z_d}(x)\to 0$ for all $x\in (0,1/2)$, as $d\to \infty$ implies uniform convergence in every set $A$ with a positive distance from points $0$ and $1$. This means that for any open set $A\subset [0,1]$ one has
$$\liminf_{d \to \infty} P(Z_d\in A)\geq  \left(\frac{1}{2} \delta_0+\frac{1}{2} \delta_1 \right)(A)$$
and the proof is complete.
\end{proof} 
\end{proposition}

The proof of Theorem \ref{ekalause} requires the following purely analytic lemma.

\begin{lemma}\label{ekalemma} Let $(g_d)_{d>0}$ be a family of increasing or decreasing functions defined on the interval $[a,b],$ where $-\infty<a<b<\infty$. Assume further that for every $d>0$ the function $g_d$ is continuously differentiable on the whole interval  $[a,b]$. Finally, assume that $0<|g_d(x)|<M$ holds for every $d>0$ and every $x\in [a,b]$.

If
\begin{equation}\label{oletus}
\lim_{d\to \infty} \left| \frac{g_d'(x)}{g_d(x)}\right|=\infty
\end{equation}
for every $x\in (a,b)$, then, for all $x\in (a,b)$,
$$g_d(x)\to0,$$
as $d\to \infty$.
\begin{proof} Without loss of generality we may assume that $(g_d)_{d>0}$ is a family of increasing and positive functions. Suppose in the contrary that there exists a number $\eta\in(a,b)$ such that
\begin{equation}\label{ehto1}
\limsup_{d\to \infty} g_d(\eta)>0. 
\end{equation}
Equation \eqref{ehto1} implies that there exists a sequence of functions $(g_{d_k})_{k=1}^\infty$, where $d_k\uparrow \infty$, as $k\to \infty$ such that $
C:=\liminf_{k\to \infty} g_{d_k}(\eta)>0$.
The fact that $g_d$ is increasing for any $d>0$ implies the inequality
\begin{equation}\label{ehto2}
\inf_{x\in [\eta,b)} \{\liminf_{k\to \infty} g_{d_k}(x)\}\geq C. 
\end{equation}
Hence, for large enough $d_k$ and all $x\in [\eta,b)$ it holds that 
\begin{equation}\label{contradiction}
\log (C/2)< \log g_{d_k}(x)<\log M.
\end{equation}

Rewriting Assumption \eqref{oletus} as
\begin{equation}\label{oletusver2}
\lim_{d\to \infty} \left| \frac{d}{dx} \log g_d(x)\right|=\infty
\end{equation}
gives
$$\lim_{k\to \infty} \frac{d}{dx} \log g_{d_k}(x)=\infty$$
for all $x\in [\eta,b)$. Set $h_{d_k}(x):=\log g_{d_k}(x)$. Now, using the fundamental theorem of calculus and the lower limit of \eqref{contradiction}, one obtains
$$h_{d_k}(b)\geq \log(C/2)+\int_{\eta}^b h_{d_k}'(y) \, dy.$$
This yields a contradiction: The function $h_{d_k}'$ is non-negative because $g_{d_k}$, and thus $\log g_{d_k}$, is increasing. Therefore, application of Fatou's lemma implies
$$\int_{\eta}^b h_{d_k}'(y) \, dy\to \infty, $$
as $k\to \infty$, contradicting the upper bound of \eqref{contradiction}.
\end{proof}
\end{lemma}

\begin{theorem}\label{ekalause} Suppose $f$ is eventually strictly log-convex or log-concave. Assume further that 
\begin{equation}\label{hyvaol}
\lim_{d \to \infty} d  \left| \frac{f'(dx)}{f(dx)}-\frac{f'(d(1-x))}{f(d(1-x))} \right|=\infty.
\end{equation}
for every $x\in(0,1/2)$.

Then $f_{Z_d}(x)\to 0$ for every $x\in (0,1/2)$.

\begin{proof} Let $x\in (0,1/2)$. Based on the proof of Lemma \ref{typechar} it is possible to choose a number $d_0$ such that the function  $f_{Z_d}$ is monotone in the interval $(x-\epsilon,x+\epsilon)\subset (0,1/2)$, when $d>d_0$, and $\epsilon>0$ is a small enough number.  

We plan to apply Lemma \ref{ekalemma} to family $(f_{Z_d})_{d>d_0}$ and interval $(a,b):=(x-\epsilon,x+\epsilon)$. To do this, note that the derivative of \eqref{f1} may be written as 
\begin{equation}\label{der1}
f_{Z_d}'(x)=d f_{Z_d}(x) \left[ \frac{f'(dx)}{f(dx)}-\frac{f'(d(1-x))}{f(d(1-x))} \right].
\end{equation}
Thus, Assumption \eqref{hyvaol} corresponds to Assumption \eqref{oletus} of Lemma \ref{ekalemma}. The remaining assumptions are clearly valid.
\end{proof}
\end{theorem}

\subsection{Main Corollary and Examples}

Theoretical results of Section \ref{tresults} imply the following surprising corollary. It is based on Proposition \ref{typechar} and Theorem \ref{ekalause}.

\begin{corollary}  There exist non-negative random variables $X$ and $Y$ such that:
\begin{enumerate} 
\item \label{part1} The variable $Y$ is asymptotically dominated by $X$, i.e. 

\begin{equation}\label{dom1}
\lim_{x \to \infty} \frac{P(Y>x)}{P(X>x)}=0,
\end{equation}
yet $Y$ is of Type \ref{a} while $X$ is of Type \ref{b}.
\item \label{part2} There exists a light-tailed random variable of Type \ref{a}.
\end{enumerate} 
\begin{proof} Define the densities $f_X$ and $f_Y$ of variables $X$ and $Y$ by 
$$f_X(x):=C_Xe^{-x+\sqrt{x}}$$
and
$$f_Y(x):=C_Ye^{-x-\sqrt{x}}$$
for $x>0$, where $C_X^{-1}=\int_0^\infty e^{-y+\sqrt{y}} \, dy$ and $C_Y^{-1}=\int_0^\infty e^{-y-\sqrt{y}} \,  dy$.

Application of L'H\^{o}pital's rule shows \eqref{dom1}. For any $x>0$,
$$\frac{d^2}{dx^2}\log f_X(x)=-\frac{1}{4} x^{-\frac{3}{2}}\, \textnormal{ and }\,  \frac{d^2}{dx^2}\log f_Y(x)=\frac{1}{4} x^{-\frac{3}{2}}. $$
Furthermore, for $x\in(0,1/2)$, we obtain
$$d  \left( \frac{f_X'(dx)}{f_X(dx)}-\frac{f_X'(d(1-x))}{f_X(d(1-x))} \right)=\frac{1}{2}\sqrt{d} ( x^{-1/2} - (1-x)^{-1/2} )\stackrel{d\to \infty}{\to} \infty$$
and
$$d  \left( \frac{f_Y'(dx)}{f_Y(dx)}-\frac{f_Y'(d(1-x))}{f_Y(d(1-x))} \right)=\frac{1}{2}\sqrt{d} (- x^{-1/2} + (1-x)^{-1/2} )\stackrel{d\to \infty}{\to} -\infty.$$

Hence, Theorem \ref{ekalause} combined with Proposition \ref{typechar} gives the result of Part \ref{part1}. The statement of Part \ref{part2} is clear because for $0<s<1$ it holds that 
$$E(e^{sY})=\int_0^\infty e^{sy}f_Y(y) \, dy<\infty $$
and thus $Y$ is a light-tailed random variable.

\end{proof}
\end{corollary}

The condition $f_{Z_d}(x)\to0$ for all $x\in (0,1/2)$ of Proposition \ref{typechar} can be difficult to verify directly. However, the sufficient condition of Theorem \ref{ekalause} seems to cover the most common situations. The class of power densities forms a notable exception. These densities are simple enough to be handled directly via Proposition \ref{typechar}. This is demonstrated in Example \ref{ex2} below.

\begin{example}\label{ex11} We check Condition \eqref{hyvaol} for certain distribution types. In all cases $x\in (0,1/2)$ and $C$ is an integration constant.
\begin{enumerate}[a)]
\item \label{weibull}Suppose $f(t)=Ce^{-t^\alpha}$, where $t>0$ and $\alpha>0$. Then
\begin{eqnarray*} && d \left( \frac{f'(dx)}{f(dx)}-\frac{f'(d(1-x))}{f(d(1-x))} \right) \\
&=&\alpha d^\alpha ((1-x)^{\alpha-1}-x^{\alpha-1})\stackrel{d\to \infty}{\to}\begin{cases} 
      \infty & \alpha>1 \\
      -\infty & 0<\alpha<1\\
      0 & \alpha=1.
   \end{cases}
\end{eqnarray*}
\item \label{lognormal} Suppose $f(t)=Ct^{-1}e^{-(\log t)^2}$, where $t>t_0>0$ for some $t_0$, $\beta\in \mathbb{R}$ and $\gamma>1$. Then
\begin{eqnarray*}
&& d  \left( \frac{f'(dx)}{f(dx)}-\frac{f'(d(1-x))}{f(d(1-x))} \right) \\
&=&2 \log d\left( \frac{1}{1-x}-\frac{1}{x}\right)+2\left( \frac{\log(1-x)+1/2}{1-x}-\frac{\log(x)+1/2}{x}\right)\stackrel{d\to \infty}{\to}-\infty.
\end{eqnarray*}
\end{enumerate}
\end{example}

Example \ref{ex11} shows that Condition \eqref{hyvaol} is satisfied by Weibull and Lognormal type densities. The next example illustrates a situation where \eqref{hyvaol} does not apply, but instead Proposition \ref{typechar} can be applied directly.
\begin{example}\label{ex2} Suppose $f(t)=t^{-\alpha}$ for some $\alpha>1$ and and all $t>t_0>0$ for some $t_0$. Then for $d>t_0/x$, where $x\in (0,1/2)$ it holds that
$$d  \left( \frac{f'(dx)}{f(dx)}-\frac{f'(d(1-x))}{f(d(1-x))} \right)=\alpha\left( \frac{1}{1-x}-\frac{1}{x}\right),$$
i.e. \eqref{hyvaol} is not valid.

However, a direct calculation using \eqref{f1} reveals that  
$$f_{Z_d}(x)=\frac{f(d x)f(d(1-x))}{\int_{0}^1 f(d y)f(d(1-y)) \, dy}\leq \frac{x^{-\alpha}(1-x)^{-\alpha}}{\int_{t_0/d}^{1-t_0/d} y^{-\alpha}(1-y)^{-\alpha} \, dy}  \to 0,$$
as $d\to \infty$.
\end{example}

\section{Discussion}\label{discussion}

Recall that $X_1$ and $X_2$ are i.i.d non-negative variables. The class of subexponential distributions $\mathcal{S}$ consists of those distributions for which
\begin{equation}\label{subexp}
\lim_{x \to \infty} \frac{P(X_1+X_2>x)}{P(X_1>x)}=2
\end{equation}
or equivalently
\begin{equation}\label{subexp2}
\lim_{x \to \infty} P(X_1>x|X_1+X_2>x)=\frac{1}{2}.
\end{equation}
In addition, the class of locally subexponential or $\Delta$-subexponential distributions $\mathcal{S}_\Delta$ can be determined by demanding that for some $\Delta>0$:
\begin{equation}\label{locsubexp}
\lim_{x \to \infty} \frac{P(X_1+X_2\in (x,x+\Delta])}{P(X_1\in (x,x+\Delta])}=2
\end{equation}
and that for any $y>0$
\begin{equation}
\lim_{x \to \infty} \frac{P(X_1\in (x,x+y+\Delta])}{P(X_1\in (x,x+\Delta])}=1.
\end{equation}
These distributions and their connections to the principle of a single big jump have been extensively studied in  \cite{watanabe:2010,teugels:1975,foss:2013,asmussen:2003,embrechts:1997,borovkov:2008}.
It is important to note that the requirement of subexponentiality or local subexponentiality does not impose detailed requirements about the distribution of $X_1$ given $X_1+X_2$. Hence, it can be argued that the process $(Z_d)$ is more suitable to describe the phenomenon of a single big jump than the membership of these distribution classes. Furthermore, it is known that $\mathcal{S}_\Delta\subset \mathcal{S}$ and that all subexponential distributions are heavy-tailed. 

In conclusion, the transition between different asymptotic Behaviours \ref{a} or \ref{b} seems to be connected to the eventual convexity or concavity of the function $\log f$. In this sense, heavy-tailedness or membership of a subexponential class has perhaps less impact on the asymptotic behaviour of $(Z_d)$ than what has been anticipated earlier.

\section*{Acknowledgements}

The deepest gratitude is expressed to the Finnish Doctoral Programme in Stochastics and Statistics (FDPSS) and the Centre of Excellence in Computational Inference (COIN) for financial support (Academy of Finland grant number 251170). Special thanks are due to Harri Nyrhinen for his diligent guidance throughout the writing of the paper.

\begin{footnotesize}

\end{footnotesize}

\end{document}